\documentclass[12pt,fullpage]{article}

\usepackage{amsfonts,xr,graphicx,amssymb,amsmath,xr,graphicx,color}

\def\R{\mathbb{R}}
\def\N{\mathbb{N}}
\def\L{\mathcal{L}}
\def\tilde{\widetilde}
\def\epsilon{\varepsilon}

\newcommand{\be}{\begin{equation}}
\newcommand{\ee}{\end{equation}}
\newcommand{\baa}{\begin{array}}
\newcommand{\eaa}{\end{array}}
\newcommand{\ba}{\begin{eqnarray}}
\newcommand{\ea}{\end{eqnarray}}

\newtheorem{theo}{\bf Theorem}[section]
\newtheorem{lem}[theo]{\bf Lemma}

\newtheorem{rem}[theo]{\bf Remark}

\begin{document}
\date{}
\title{\bf{Convexity of level sets for elliptic problems in convex domains or convex rings: two counterexamples}}
\author{Fran{\c{c}}ois Hamel$^{\hbox{\small{ a,b}}}$, Nikolai Nadirashvili$^{\hbox{\small{ a}}}$ and Yannick Sire$^{\hbox{\small{ a}}}\,$\thanks{The research leading to these results has received funding from the French ANR within the project PRE\-FERED and from the European Research Council under the European Union's Seventh Framework Programme (FP/2007-2013) / ERC Grant Agreement n.321186 - ReaDi - Reaction-Diffusion Equations, Propagation and Modelling. Part of this work was also carried out during visits by the first author to the Departments of Mathematics of the University of California, Berkeley and of Stanford University, the hospitality of which is thankfully acknowledged.}\\
\\
\footnotesize{$^{\hbox{a }}$Aix-Marseille Universit\'e}\\
\footnotesize{LATP (UMR CNRS 7353), 39 rue F. Joliot-Curie, F-13453 Marseille Cedex 13, France}\\
\footnotesize{$^{\hbox{b }}$Institut Universitaire de France}}
\maketitle

\begin{abstract}
\noindent{}This paper deals with some geometrical properties of solutions of some semilinear elliptic equations in bounded convex domains or convex rings. Constant boundary conditions are imposed on the single component of the boundary when the domain is convex, or on each of the two components of the boundary when the domain is a convex ring. A function is called quasiconcave if its superlevel sets, defined in a suitable way when the domain is a convex ring, are all convex. In this paper, we prove that the superlevel sets of the solutions do not always inherit the convexity or ring-convexity of the domain. Namely, we give two counterexamples to this quasiconcavity property: the first one for some two-dimensional convex domains and the second one for some convex rings in any dimension.
\end{abstract}

%%%%%%%%%%%%%%%%%%%%%%%%%%%%%%%%%%%%%%%%%
%%%%%%%%%%%%%%%%%%%%%%%%%%%%%%%%%%%%%%%%%

\section{Introduction and main results}\label{intro}

This paper is concerned with some geometrical properties of real-valued solutions of semilinear elliptic equations
\be\label{eq}
\Delta u+f(u)=0
\ee
in bounded domains $\Omega\subset\R^N$, in dimensions $N=2$ or higher, with Dirichlet-type boundary conditions on~$\partial\Omega$. By domains, we mean non-empty open connected subsets of $\R^N$.

The domains~$\Omega$ are assumed to be either convex domains or convex rings. One is interested in knowing how these geometrical properties of~$\Omega$ are inherited by the solutions~$u$, under some suitable boundary conditions, that is how the shape of the solutions is influenced by the shape of the underlying domains.  It is well-known that the convexity or the concavity of the solutions are too strong properties which are not true in general (see e.g.~\cite{ko1}). However, a typical question we address in this paper is the following one: assuming that~$\Omega$ is convex and that $u$ is a solution of~(\ref{eq}) which is positive in~$\Omega$ and vanishes on $\partial\Omega$, is it true that the superlevel sets
$$\big\{x\in\Omega;\ u(x)>\lambda\big\}$$
of~$u$ are all convex? A similar question can be asked when $\Omega$ is a convex ring and~$u$ is assumed to be equal to two constant values on the two connected components of $\partial\Omega$ (see below for detailed statements). These questions have been well studied in the literature: more precisely, almost all papers on this field have been devoted to the proof of a positive answer to these questions, under some suitable conditions on the function $f$, see the references below. In this paper, we prove that the answer to these questions can also be negative, that is we show that the superlevel sets of some solutions $u$ of problems of the type~(\ref{eq}) are not all convex. More precisely, we give two counterexamples, one in a class of convex domains and one in a class of convex rings.

Let us first deal with the case of bounded convex domains~$\Omega$ and let us consider the semilinear elliptic problem
\be\label{pbconvex}\left\{\baa{rcll}
\Delta u+f(u) & = & 0 & \hbox{in }\Omega,\vspace{3pt}\\
u & = & 0 & \hbox{on }\partial\Omega,\vspace{3pt}\\
u & > & 0 & \hbox{in }\Omega.\eaa\right.
\ee
Throughout the paper, the function $f:[0,+\infty)\to\R$ is assumed to be locally H\"older continuous. The domains~$\Omega$ are always assumed to be of class~$C^{2,\alpha}$ (with $\alpha>0$, we then say that the domains~$\Omega$ are smooth) and the solutions $u$ are understood in the classical sense~$C^2(\overline{\Omega})$. The superlevel set~$\big\{x\in\Omega;\ u(x)>0\big\}$ of a solution $u$ of~(\ref{pbconvex}) is equal to the domain $\Omega$, which is convex by assumption. A natural question is to know whether the superlevel sets $\big\{x\in\Omega;\ u(x)>\lambda\big\}$ for~$\lambda\ge 0$ are all convex or not. If this is the case,~$u$ is called quasiconcave.

In his paper~\cite{li} (see Remark~3, page~268), P.-L.~Lions writes that, in a convex domain~$\Omega$, ``[he] believe[s] that [...] for general $f$, the [super]level sets of any solution $u$ of~[(\ref{pbconvex})] are convex". There is indeed a vast literature containing some proofs of the above statement for various nonlinearities~$f$. We here list some of the most classical references. Firstly, Makar-Limanov~\cite{ml} proved that, for the two-dimensional torsion problem, that is $f(u)=1$ with $N=2$, the solution $u$ is quasiconcave, since~$\sqrt{u}$ is actually concave. Brascamp and Lieb~\cite{bl} showed that, if~$f(u)=\lambda u$ ($\lambda$ is then necessarily the principal eigenvalue of the Laplacian with Dirichlet boundary condition), then the principal eigenfunction $u$ is quasiconcave and more precisely it is log-concave, that is $\log u$ is concave. The proof uses the fact that log-concavity is preserved by the heat equation (but quasiconcavity is not in general, see~\cite{is}). When~$f(u)=\lambda u^p$ with $0<p<1$ and $\lambda>0$, Keady~\cite{k} for $N=2$ and Kennington~\cite{ke1} for $N\ge 2$ proved that $u^{(1-p)/2}$ is concave, whence $u$ is quasiconcave. Many generalizations under more general assumptions on $f$ and alternate proofs have been given. A possible strategy is to prove that $g(u)$ is concave for some suitable increasing function $g$, by showing that
$$g(u(tx+(1-t)y))-tg(u(x))-(1-t)g(u(y))\ge0\hbox{ for all }(t,x,y)\in[0,1]\times\Omega\times\Omega$$
and by using the elliptic maximum principle or the preservation of concavity of $g(u)$ by a suitable parabolic equation, see~\cite{casp,gp,ka3,ka4,ke1,ke2,ko1,li}. Other strategies consist in studying the sign of the curvatures of the level sets of $u$ or in proving that the Hessian matrix of $g(u)$ for some suitable increasing $g$ has a constant rank, see~\cite{app,bgmx,cf,kl,l,x}. Lastly, we refer to~\cite{all,gg} for further references using the quasiconcave envelope and singular perturbations arguments, and to the book of Kawohl~\cite{ka2} for a general overview.

The first main result of this paper is, to our best knowledge, the first counterexample to the quasiconcavity of solutions $u$ of~(\ref{pbconvex}) in convex domains $\Omega$.

\begin{theo}\label{th1} In dimension $N=2$, there are some smooth bounded convex domains~$\Omega$ and some~$C^{\infty}$ functions $f:[0,+\infty)\to\R$ such that
$$f(s)\ge 1\ \hbox{ for all }s\ge 0$$
and for which problem~$(\ref{pbconvex})$ admits both a quasiconcave solution~$v$ and a solution~$u$ which is not quasiconcave.
\end{theo}

\begin{rem}{\rm When $\Omega$ is an Euclidean ball of $\R^N$ in any dimension $N\ge 1$ and when $f$ is locally Lipschitz-continuous, then the celebrated paper of Gidas, Ni and Nirenberg~\cite{gnn} asserts that any solution~$u$ of~$(\ref{pbconvex})$ is radially symmetric and decreasing with respect to the center of the ball: in other words, the superlevel sets of~$u$ are all concentric balls and $u$ is therefore quasiconcave. In particular, Theorem~\ref{th1} cannot hold in dimension $N=1$. More generally speaking, if the convex domain~$\Omega$ is symmetric with respect to some hyperplane, then the moving plane method implies that~$u$ itself inherits this property and is actually symmetric and decreasing with respect to the distance to this hyperplane, see~\cite{gnn}. As a matter of fact, the two-dimensional convex domains~$\Omega$ constructed in the proof of Theorem~\ref{th1} are symmetric with respect to both variables~$x$ and~$y$ of $\R^2$, whence any solution $u$ of~$(\ref{pbconvex})$ is symmetric with respect to~$x$ and~$y$, and decreasing with respect to $|x|$ and $|y|$. These properties imply that the superlevel sets of~$u$ are necessarily symmetric and convex with respect to~$x$ and~$y$, and starshaped with respect to the origin $(0,0)$. But the symmetry and convexity properties of the superlevel sets in $x$ and $y$ do not mean that these superlevel sets are truly convex! Actually, they are not so in general, as Theorem~$\ref{th1}$ shows.}
\end{rem}

\begin{rem}\label{remstability}{\rm In~\cite{cc}, Cabr\'e and Chanillo proved that, if $\Omega$ is a smooth bounded strictly convex domain of $\R^2$ and if $u$ is any semi-stable solution of~$(\ref{pbconvex})$ in the sense that
$$\int_{\Omega}|\nabla\phi|^2-\int_{\Omega}f'(u)\phi^2\ge0$$
for every $C^{\infty}(\Omega)$ function $\phi$ whose support is compactly included in $\Omega$, then $u$ has a unique critical point $($its maximum$)$ and this critical point is nondegenerate, whence the superlevel sets $\Omega^{\lambda}$ of~$u$ are convex for $\lambda$ close to $\max_{\overline{\Omega}}u$ (for further results about the uniqueness and nondegeneracy of the critical point in some more general convex domains~$\Omega$, we refer to~\cite{ag,cc,p,sp}). If the semi-stability were known to imply the convexity of {\it all} superlevel sets $($and also in convex domains which are not {\it strictly} convex$)$, then the solutions~$u$ constructed in Theorem~$\ref{th1}$ would therefore not be semi-stable. However, proving the quasiconcavity from the semi-stability is still an open question, as well as proving or disproving directly the semi-stability of the solutions $u$ of Theorem~$\ref{th1}$. Notice that if $f'(u)$ were nonpositive in $\Omega$, then $u$ would be automatically semi-stable. For the solutions~$u$ of Theorem~$\ref{th1}$, the function $f'(u)$ is actually equal to $0$ on a large set, and the set $\big\{x\in\Omega;\ f'(u(x))>0\big\}$ is always a non-empty open set $($but one can not directly infer the unstability of $u$ from this sole property$)$. Lastly, the set $\big\{x\in\Omega;\ f'(u(x))<0\big\}$ is not empty in general for the solutions~$u$ of Theorem~$\ref{th1}$, whence $f'$ has in general no sign and $f$ is in general non-monotone on the range~$[\min_{\overline{\Omega}}u,\max_{\overline{\Omega}}u]=[0,\max_{\overline{\Omega}}u]$ of $u$ $($see the proof of Theorem~$\ref{th1}$ and Remark~$\ref{remf'}$ for more details$)$.}
\end{rem}

In the second part of the paper, we deal with the case of convex rings~$\Omega$ in any dimension~$N\ge 2$. Namely, a domain $\Omega\subset\R^N$ is called a convex ring if
$$\Omega=\Omega_1\backslash\overline{\Omega_2},$$
where $\Omega_1$ and $\Omega_2$ are two bounded convex domains of $\R^N$ such that
$$\overline{\Omega_2}\subset\Omega_1.$$
In a convex ring~$\Omega=\Omega_1\backslash\overline{\Omega_2}$, let us now consider the semilinear elliptic problem
\be\label{pbring}\left\{\baa{rcll}
\Delta u+f(u) & \!\!=\!\! & 0 & \hbox{in }\Omega,\vspace{3pt}\\
u & \!\!=\!\! & 0 & \hbox{on }\partial\Omega_1,\vspace{3pt}\\
u & \!\!=\!\! & M & \hbox{on }\partial\Omega_2,\vspace{3pt}\\
u & \!\!>\!\! & 0 & \hbox{in }\Omega,\eaa\right.
\ee
where $M>0$ is a positive real number. For any classical solution $u$ of~(\ref{pbring}), we define the function~$\overline{u}\in C(\overline{\Omega_1})$ by
$$\overline{u}(x)=\left\{\baa{ll}
u(x) & \hbox{if }x\in\overline{\Omega},\\
M & \hbox{if }x\in\Omega_2=\overline{\Omega_1}\backslash\overline{\Omega},\eaa\right.$$
and we say that $u$ is quasiconcave in $\Omega$ if $\overline{u}$ is so in $\Omega_1$, that is if the superlevel sets
$$\Omega^{\lambda}:=\big\{x\in\Omega_1;\ \overline{u}(x)>\lambda\big\}$$
are convex for all $\lambda\ge 0$. Notice that $\Omega^0=\Omega_1$ is convex by assumption. If we knew that~$u<M$ in $\Omega$, then $\bigcap_{\lambda<M}\Omega^{\lambda}=\big\{x\in\Omega_1;\ \overline{u}(x)\ge M\big\}=\overline{\Omega_2}$ would be convex too. However, the condition~$u<M$ in $\Omega$ is not imposed a priori, even if the solutions $u$ which will be constructed in Theorem~\ref{th2} below satisfy this property under some specific conditions on $f$ (see Remarks~\ref{remmu} and~\ref{remM}).

Many papers have been devoted to finding sufficient conditions on $f$ which guarantee the convexity of the superlevel sets of the solutions $u$ of~(\ref{pbring}). The first positive classical result of Gabriel~\cite{g} is concerned with three-dimensional harmonic functions (see also Lewis~\cite{le} for $p$-harmonic functions). Caffarelli and Spruck~\cite{casp} proved the quasiconcavity of $u$ in any dimension~$N\ge 2$ when $f(0)=0$ and $f$ is nonincreasing. We refer to~\cite{cf,ko2,kl} for further positive results using properties of the curvatures of the level sets of $u$ or the rank of the Hessian matrix of $g(u)$ for some increasing $g$, to~\cite{casp,dk,gr,ka1} for positive results using properties of minimal points of the quasiconcavity function $(x,y)\mapsto u((x+y)/2)-\min(u(x),u(y))$ in $\Omega\times\Omega$, to~\cite{bls,cosa,cusa} for positive results using the maximum principle for the quasiconcave envelope of the function~$u$ and to~\cite{aps,ls} for further existence results of quasiconcave solutions to some equations of the type~(\ref{pbring}). Lastly, if the open sets $\Omega_1$ and $\Omega_2$ are just assumed to be starshaped with respect to a point~$x_0\in\Omega_2$, then the superlevel sets of the solutions~$u$ of~(\ref{pbring}) are known to be starshaped with respect to~$x_0$ when~$f(0)=0$ and $f$ is nonincreasing, since $(x-x_0)\cdot\nabla u(x)<0$ in $\Omega$ from the maximum principle and Hopf lemma, see~\cite{a2,dk,gm,ka1,sa} for further results in this direction.

The only counterexample to the quasiconcavity of solutions $u$ of~(\ref{pbring}), to our best knoweldge, is the one of Monneau and Shahgholian~\cite{ms}: the authors prove that, in dimension $N=2$, for some convex rings and for some nonnegative functions $f$ which are close to a Dirac mass concentrated at some real number between~$0$ and $M$, the solutions of~(\ref{pbring}) in $\Omega$ cannot be quasiconcave. The construction uses the existence of non-convex domains solving some approximated free boundary problems (see~\cite{a1}).

The second main result of the present paper gives a counterexample to the quasiconcavity of solutions $u$ of~(\ref{pbring}), of a type different from~\cite{ms}. The function~$f:[0,+\infty)\to\R$ will be any locally Lipschitz-continuous function such that
\be\label{hypf}\left\{\baa{l}
f\hbox{ is bounded from above, that is }\displaystyle{\mathop{\sup}_{s\in[0,+\infty)}}\,f(s)<+\infty,\vspace{3pt}\\
s\mapsto\displaystyle{\frac{f(s)}{s}}\hbox{ is decreasing over }(0,+\infty),\vspace{3pt}\\
\hbox{either }f(0)>0,\ \hbox{ or }\ f(0)=0\hbox{ and }\displaystyle{\mathop{\lim}_{s\to 0^+}}\,\displaystyle{\frac{f(s)}{s}}>\lambda_1(-\Delta,\Omega_1)\,(>0),\eaa\right.
\ee
where $\lambda_1(-\Delta,\Omega_1)$ denotes the smallest eigenvalue of the operator $-\Delta$ in $\Omega_1$ with Dirichlet boun\-dary condition on $\partial\Omega_1$.

\begin{theo}\label{th2}
Let $N$ be any integer such that $N\ge 2$, let $\Omega_1$ be any smooth bounded convex domain of~$\R^N$ and let $f$ be any function satisfying~$(\ref{hypf})$. Then there exists a constant $M_0>0$ such that, for all $M\ge M_0$, there are some smooth convex rings $\Omega=\Omega_1\backslash\overline{\Omega_2}$ for which problem~$(\ref{pbring})$ has a unique solution~$u$ and this solution $u$ is not quasiconcave.
\end{theo}

In Theorem~\ref{th2}, the domain $\Omega_1$ is any given convex domain and $f$ is any given function satisfying~(\ref{hypf}). One of the main assumptions, quite different from the construction given in~\cite{ms}, is that~$f'(0^+)$ is not too small if $f(0)=0$. However, even if the function~$f$ is assumed to be positive in a right-neighborhood of $0$, it may not be nonnegative everywhere. Typical examples of functions~$f$ satisfying~(\ref{hypf}) are the positive constants $f(s)=\beta>0$, or functions of the type
\be\label{fgamma}
f(s)=\gamma\,s-s^p
\ee
with $p>1$ and $\gamma>\lambda_1(-\Delta,\Omega_1)$ (when $\gamma$ is a fixed positive constant, this last condition is automatically fulfilled if $\Omega_1$ contains a ball with a large enough radius). Notice that for nonlinearities of the type~(\ref{fgamma}) with $\gamma>\lambda_1(-\Delta,\Omega_1)$, there exists a unique solution $u$ of problem~(\ref{pbconvex}) in~$\Omega_1$ (see~\cite{b}) and this solution is log-concave, whence quasiconcave (see~\cite{li}).

We also point out that Theorem~\ref{th2} holds in any dimension. As a matter of fact, it also holds for equations which are much more general than~(\ref{pbring}), with non-symmetric operators or heterogeneous coefficients. For the sake of clarity of the presentation we prefered to state only the counterexamples for problem~(\ref{pbring}) in the present section. We refer to Section~\ref{secrings} for more general problems.

\begin{rem}{\rm In problem~$(\ref{pbring})$ and in Theorem~$\ref{th2}$, one can replace the boundary condition $u=M$ on $\partial\Omega_2$ by $u=1$ $($or by any other arbitrary positive real number$)$, even if it means changing $f$. More precisely, if~$f$, $\Omega_1$ and $M_0>0$ are as in Theorem~$\ref{th2}$, then, for every $M\ge M_0$, the function~$\tilde{u}=u/M$ indeed solves
$$\left\{\baa{rcll}
\Delta\tilde{u}+\tilde{f}(\tilde{u}) & \!\!=\!\! & 0 & \hbox{in }\Omega,\vspace{3pt}\\
\tilde{u} & \!\!=\!\! & 0 & \hbox{on }\partial\Omega_1,\vspace{3pt}\\
\tilde{u} & \!\!=\!\! & 1 & \hbox{on }\partial\Omega_2,\vspace{3pt}\\
\tilde{u} & \!\!>\!\! & 0 & \hbox{in }\Omega,\eaa\right.$$
where $u$ and $\Omega_2$ as in the statement of Theorem~$\ref{th2}$ and $\tilde{f}(s)=f(Ms)/M$ satisfies the same condition~$(\ref{hypf})$ as $f$.}
\end{rem}

\begin{rem}\label{remmu}{\rm If, in addition to~$(\ref{hypf})$, the function $f$ is assumed to be nonpositive for large $s$, that is there exists a real number $\mu>0$ such that $f(s)\le 0$ for all $s\ge\mu$, then one can take $M_0=\mu$ in Theorem~$\ref{th2}$ and the solutions $u$ of~$(\ref{pbring})$ are such that
$$0<u<M\ \hbox{ in }\Omega$$
for all $M\ge\mu$. We refer to Remark~$\ref{remM}$ for further details.}
\end{rem}

%%%%%%%%%%%%%%%%%%%%%%%%%%%%%%%%%%%%%%%%%
%%%%%%%%%%%%%%%%%%%%%%%%%%%%%%%%%%%%%%%%%

\section{Counterexamples in convex domains: proof of Theorem~\ref{th1}}\label{secdomains}

This section is devoted to the proof of Theorem~\ref{th1}. That is, we construct explicit examples of bounded smooth convex two-dimensional domains $\Omega$ and of functions $f$ for which problem~(\ref{pbconvex}) admits some non-quasiconcave solutions~$u$. The construction is divided into five main steps. Firstly, we define a one-parameter family~$(\Omega_a)_{a\ge 1}$ of more and more elongated stadium-like convex domains. Secondly, for each value of the parameter $a\ge 1$, we solve a  variational problem in $H^1_0(\Omega_a)$ with a nonlinear constraint, any solution~$u_a$ of which solves an elliptic equation of the type~(\ref{pbconvex}) in $\Omega_a$ with some function~$f_a$. Thirdly, we prove some a priori estimates for the superlevel sets of the functions $u_a$. Next, we compare $u_a$ with a one-dimensional profile in $\Omega_a$ when $a$ is large enough. Lastly, we show that the superlevel sets of the functions $u_a$ cannot be all convex when $a$ is large enough.\par
As a preliminary step, let us fix a $C^{\infty}$ function $g:\R\to[0,1]$ such that
\be\label{choiceg}
g=0\hbox{ on }(-\infty,1],\ \ g=1\hbox{ on }[2,+\infty)\ \hbox{ and }\ g'\ge0\hbox{ on }\R.
\ee
The function $g$ is fixed throughout the proof.\hfill\break\par

{\it Step 1: construction of a family of smooth bounded convex domains~$(\Omega_a)_{a\ge 1}$.} We first introduce a family of stadium-like smooth convex domains. Denote $(x,y)$ the coordinates in $\R^2$. Let~$\varphi:[-1,1]\to\R$ be a fixed continuous nonnegative concave even function such that $\varphi(\pm 1)=0$. For $a\ge 1$, we define
\be\label{defOmegaa}
\Omega_a=\big\{(x,y)\in\R^2;\ -a-\varphi(y)<x<a+\varphi(y),\ -1<y<1\big\}
\ee
and we choose $\varphi$ once for all so that $\Omega_1$ (and then $\Omega_a$ for every $a\ge 1$) be of class $C^{2,\alpha}$ with~$\alpha>0$ (this means that $\varphi$ is of class $C^{2,\alpha}_{loc}(-1,1)$ and that $\varphi$ satisfies some compatibility conditions at~$\pm 1$). The $C^{2,\alpha}$ bounded domains $\Omega_a$ for $a\ge 1$ are all convex and axisymmetric with respect to both axes~$\{x=0\}$ and $\{y=0\}$, see the joint figure.\par
\begin{figure}\label{figure1}
\begin{center}
\includegraphics[scale=0.85]{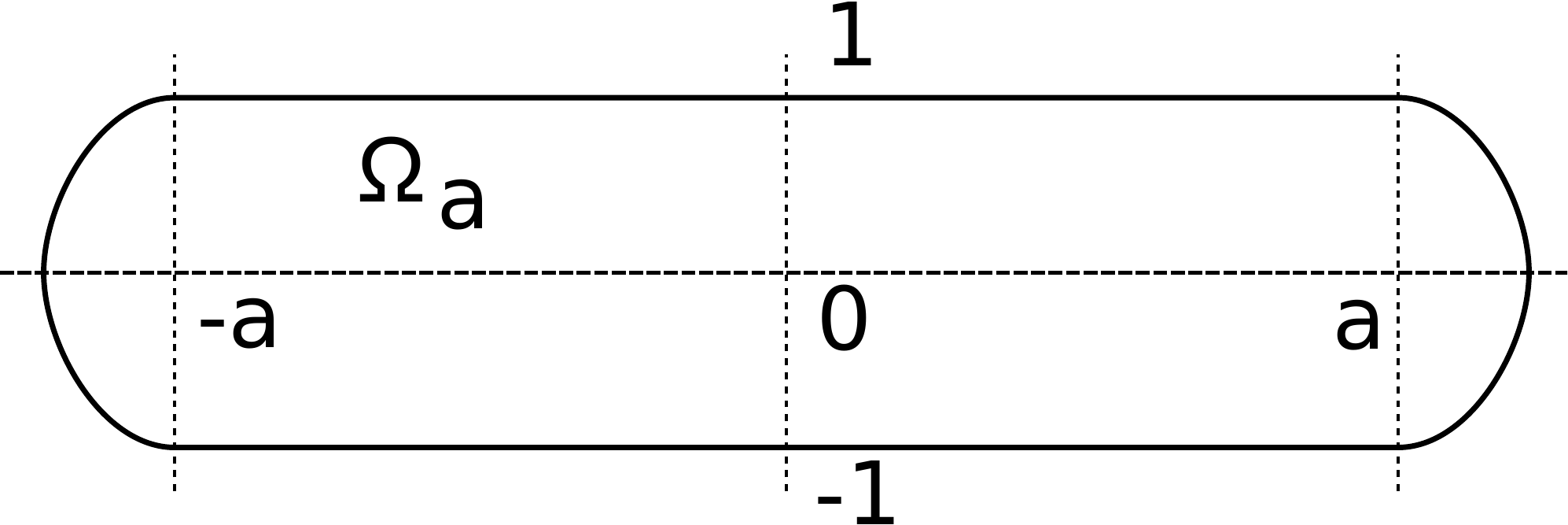}
\caption{The convex stadium-like domain $\Omega_a$}
\end{center}
\end{figure}
Our goal is to show that the conclusion of Theorem~\ref{th1} holds with these convex domains $\Omega_a$ and some functions $f_a$, when $a$ is large enough.\hfill\break\par

{\it Step 2: a constrained variational problem in $\Omega_a$.} In this step, we fix a parameter $a\ge 1$. We construct a $C^{2,\alpha}(\overline{\Omega_a})$ function $u_a$ as a minimizer of a constrained variational problem in~$\Omega_a$.\par
Let $I_a$ be the functional defined in $H^1_0(\Omega_a)$ by
$$I_a(u)=\frac{1}{2}\int_{\Omega_a}|\nabla u|^2-\int_{\Omega_a}u,\ \ u\in H^1_0(\Omega_a).$$
It is well-known that this functional has a unique minimizer in $H^1_0(\Omega_a)$, which is the classical~$C^{2,\alpha}(\overline{\Omega_a})$ solution~$v_a$ of the torsion problem
\be\label{eqva}\left\{\baa{rcl}
\Delta v_a+1 & = & 0\ \mbox{ in }\Omega_a,\vspace{3pt}\\
v_a & = & 0\ \mbox{ on }\partial\Omega_a.\eaa\right.
\ee
It follows from the strong maximum principle and the definition of $\Omega_a$ that
\be\label{ineqva}
0<v_a(x,y)<\frac{1-y^2}{2}\ \hbox{ for all }(x,y)\in\Omega_a.
\ee
This function $v_a$ is also known to be quasiconcave in $\Omega_a$, see~\cite{ml}.\par
We are then going to replace $v_a$ by a function $u_a$ which minimizes the functional~$I_a$ over a nonlinear subset of~$H^1_0(\Omega_a)$ and which will be our non-quasiconcave candidate for a problem of the type~(\ref{pbconvex}).\par
To do so, let us now define
$$U_a=\Big\{u\in H^1_0(\Omega_a);\ \int_{\Omega_a}g(u)=1\Big\}.$$
Since the Lebesgue measure $|\Omega_a|$ of $\Omega_a$ is larger than $1$, the set $U_a$ is not empty: for instance, by continuity of the map $\R\ni t\mapsto\int_{\Omega_a}g(tv_a)$, there is a real number $t_a\in(0,+\infty)$ such that $t_av_a\in U_a$. Furthermore, it is straightforward to check, using Poincar\'e's inequality together with Rellich's and Lebesgue's theorems, that the minimum of the functional~$I_a$ over the set~$U_a$ is reached, by a function~$u_a\in U_a$, that is
$$I_a(u_a)=\min_{u\in U_a}I_a(u).$$\par
We observe that $g'(u_a)\in L^{\infty}(\Omega_a)$ is not the zero function. Otherwise, the gradient $\nabla g(u_a)$ of the $H^1(\Omega_a)$ function~$g(u_a)$ would be equal to $\nabla g(u_a)=g'(u_a)\nabla u_a=0$ a.e. in $\Omega_a$ and, by definition of~$U_a$,~$g(u_a)$ would then be equal to the positive constant $1/|\Omega_a|$ a.e. in $\Omega_a$. Due to~(\ref{choiceg}), there would then exist $m>0$ such that $u_a\ge m$ a.e. in $\Omega_a$, contradicting the fact that $u_a\in H^1_0(\Omega_a)$ has a zero trace on $\partial\Omega_a$. Hence, $g'(u_a)$ cannot be the zero function and the differential of the map~$H^1_0(\Omega_a)\ni u\mapsto\int_{\Omega_a}g(u)$ is not zero at $u_a$.\par
From the Euler-Lagrange formulation and elliptic regularity theory, any such minimizer $u_a$ is then a classical~$C^{2,\alpha}(\overline{\Omega_a})$ solution of  an equation of the type
\be\label{equa}\left\{\baa{rcl}
\Delta u_a+f_a(u_a) & = & 0\ \hbox{ in }\Omega_a,\vspace{3pt}\\
u_a & = & 0\ \hbox{ on }\partial\Omega_a,\eaa\right.
\ee
where
$$f_a(s)=1+\mu_a\,g'(s)\ \hbox{ for }s\in\R$$
and $\mu_a\in\R$ is a Lagrange multiplier. Observe that the function $f_a$ is of class~$C^{\infty}(\R)$. Furthermore,
$$\Delta(u_a-v_a)=-\mu_ag'(u_a)$$
has a constant sign in~$\Omega_a$, since~$g'$ is nonnegative. As a consequence of the maximum principle, the function $u_a-v_a$ itself has a constant sign in $\Omega_a$. But
\be\label{maxua}
\max_{\overline{\Omega_a}}u_a>1
\ee
because of~(\ref{choiceg}) and by definition of $U_a$. Therefore, from~(\ref{ineqva}), the function $v_a$ cannot majorize~$u_a$. The strong maximum principle finally implies that
\be\label{vaua}
0<v_a(x,y)<u_a(x,y)\ \hbox{ for all }(x,y)\in\Omega_a.
\ee
Thus, the function $u_a$ is a classical solution of the problem~(\ref{pbconvex}) in $\Omega_a$ with the function~$f_a$. Notice also that the sign of $\Delta(u_a-v_a)$ is therefore nonpositive and, since $u_a$ and $v_a$ are not identically equal, one has $\mu_a>0$. In particular,
\be\label{fa}
f_a(s)\ge1\ \hbox{ for all }s\in\R.
\ee
On the other hand, since $f_a(s)=1$ for all $s\ge 2$ because of~(\ref{choiceg}), the maximum principle also yields
\be\label{uasup}
u_a(x,y)<\frac{1-y^2}{2}+2\ \hbox{ for all }(x,y)\in\Omega_a.
\ee\par
The uniqueness of the minimizer $u_a$ of $I_a$ in the set $U_a$ is not clear, and is anyway not needed in the sequel. However, we point out an important geometrical property fulfilled by $u_a$, which will be used in the next step. Namely, since~$\Omega_a$ is convex and symmetric with respect to the axes $\{x=0\}$ and $\{y=0\}$, it follows from~\cite{gnn} that $u_a$ is even in $x$ and $y$ and is decreasing with respect to~$|x|$ and~$|y|$.\par
In the sequel, we are going to show that, for $a$ large enough, the conclusion of Theorem~\ref{th1} holds with $\Omega_a$, $f_a$ and $u_a$, that is the minimizers~$u_a$ have some non-convex superlevel sets. Notice that $f_a$ satisfies~(\ref{fa}), as stated in Theorem~\ref{th1}.\par
Before going further on, we also point out that the solution $v_a$ of the torsion problem~$(\ref{eqva})$ also solves the same equation~$(\ref{pbconvex})$ as $u_a$, with $f_a$ in $\Omega_a$, because of~$(\ref{ineqva})$ and the fact that $f_a=1$ on $[0,1]\supset[0,1/2]$ due to~$(\ref{choiceg})$. Therefore, problem~$(\ref{pbconvex})$ with $f_a$ in $\Omega_a$ admits the solution $v_a$, which is always quasiconcave by~\cite{ml} applied to~$(\ref{eqva})$, whereas the solutions $u_a$ will be proved to be non-quasiconcave for $a$ large.\hfill\break\par

{\it Step 3: a priori estimates of the size of a superlevel set of the functions $u_a$.} In this step, we study the location of the superlevel sets
$$\omega_a=\big\{(x,y)\in\Omega_a;\ u_a(x,y)>1\big\}$$
of the minimizers $u_a$ of $I_a$ in $U_a$ when $a$ is large. From~(\ref{maxua}) and the remarks of the previous step, the sets~$\omega_a$ are non-empty open sets, they are all symmetric with respect to the axes $\{x=0\}$ and~$\{y=0\}$, and they are convex with respect to both variables $x$ and $y$.\par
The key-point in this step is to show a uniform control of the size of the sets $\omega_a$. We first begin with a bound in the $x$-direction, meaning that the sets $\omega_a$ are not too elongated.

\begin{lem}\label{lemtaille1}
There exists a constant $C_x>0$ such that 
\be\label{taille1}
0\le\sup_{(x,y)\in\omega_a}|x|<C_x
\ee
for all $a\ge 1$ and for any minimizer $u_a$ of $I_a$ in $U_a$.
\end{lem}

\begin{figure}\label{figure2}
\begin{center}
\includegraphics[scale=0.85]{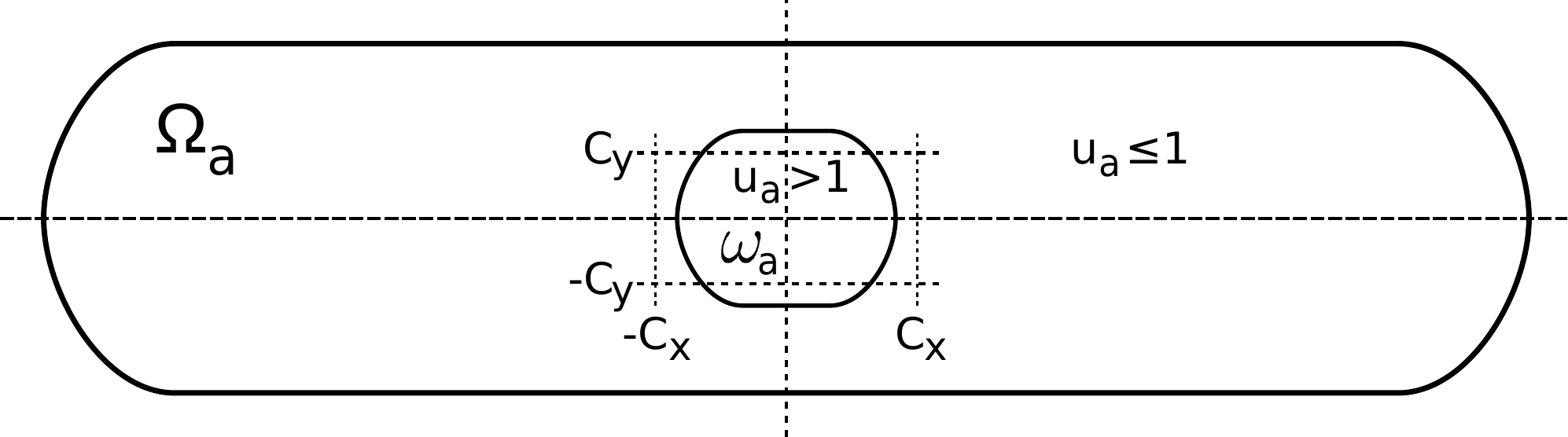}
\caption{The set $\omega_a$ where $u_a>1$}
\end{center}
\end{figure}

\noindent{\bf{Proof.}} The proof is divided into two main steps. We first estimate from above the quantities~$I_a(u_a)$ by introducing a suitable test function in the set~$U_a$, which is not too far from the one-dimensional function $y\mapsto(1-y^2)/2$. Then, we estimate $I_a(u_a)$ from below by observing that if~$u_a(x,0)$ is larger than~$1$ then the contribution of $u_a(x,\cdot)$ to $I_a(u_a)$ in the section $\Omega_a\cap(\{x\}\times\R)$ will be uniformly larger than that of the minimizer $y\mapsto(1-y^2)/2$. This eventually provides a control of the size of such points $x$ and then of the size of $\omega_a$, independently of $a$.\par
Throughout the proof, one can assume without loss of generality that $a$ is any real number such that $a\ge 2$ (since $\sup_{(x,y)\in\Omega_a}|x|\le a+\|\varphi\|_{L^{\infty}(-1,1)}$ for all $a\ge 1$ by the definition~(\ref{defOmegaa}) of $\Omega_a$). We consider any minimizer $u_a$ of the functional $I_a$ in the set $U_a$ and we set
\be\label{defxa}
x_a=\sup_{(x,y)\in\omega_a}|x|.
\ee\par
Let us first bound $I_a(u_a)$ from above by using the minimality of $u_a$ and comparing $I_a(u_a)$ with the value of $I_a$ at some suitably chosen test function. Let $w$ be a fixed $C^{\infty}(\R^2)$ nonnegative function such that
$$w=0\hbox{ in }\R^2\backslash(-1,1)^2\ \hbox{ and }\ w>0\hbox{ in }[-2/3,2/3]^2.$$
The function $w$ is independent of $a$. Let $\phi_0$ be the $H^1_0(-1,1)$ function defined by
\be\label{defphi0}
\phi_0(y)=\frac{1-y^2}{2}\ \hbox{ for all }y\in[-1,1].
\ee
We point out that $\phi_0$ is the unique minimizer in $H^1_0(-1,1)$ of the functional $J$ is defined by
\be\label{defJ}
J(\phi)=\frac{1}{2}\int_{-1}^1\phi'(y)^2dy-\int_{-1}^1\phi(y)dy,\ \ \phi\in H^1_0(-1,1).
\ee
From Lebesgue's dominated convergence theorem, the function
$$G:t\mapsto\int_{(-1,1)^2}g(\phi_0(y)+t\,w(x,y))\,dx\,dy$$
is continuous in $\R$. Furthermore, $G(0)=0$ from~(\ref{choiceg}) and~(\ref{defphi0}), and
$$\lim_{t\to+\infty}G(t)=\int_{\{w(x,y)>0\}}dx\,dy\ge\Big(\frac{4}{3}\Big)^2>1.$$
Therefore, there is $t_0\in(0,+\infty)$, independent of $a$, such that
$$G(t_0)=\int_{(-1,1)^2}g(\phi_0(y)+t_0w(x,y))\,dx\,dy=1.$$\par
Let us now consider the test function $w_a$ defined in $\overline{\Omega_a}$ by
$$w_a(x,y)=\phi_0(y)\chi_a(x)+t_0w(x,y),$$
where $\chi_a:\R\to[0,1]$ is even and defined in $[0,+\infty)$ by
$$\chi_a(x)=\left\{\baa{ll}
1 & \hbox{if }x\in[0,a-1],\vspace{3pt}\\
a-x & \hbox{if }x\in(a-1,a),\vspace{3pt}\\
0 & \hbox{if }x\ge a.\eaa\right.$$
The function $w_a$ belongs to $H^1_0(\Omega_a)$. Furthermore, since $a\ge 2$, one has
$$w_a(x,y)=\phi_0(y)+t_0w(x,y)\ \hbox{ for all }(x,y)\in(-1,1)^2,$$
while
$$w_a(x,y)=\phi_0(y)\chi_a(x)\le\phi_0(y)<1\ \hbox{ for all }(x,y)\in\Omega_a\backslash(-1,1)^2.$$
Therefore,
$$\int_{\Omega_a}g(w_a)=\int_{(-1,1)^2}g(w_a)=\int_{(-1,1)^2}g(\phi_0(y)+t_0w(x,y))\,dx\,dy=G(t_0)=1.$$
In other words, $w_a\in U_a$. By definition of $u_a$, one infers that
\be\label{Ianunvn}
I_a(u_a)\le I_a(w_a).
\ee\par
Let us now estimate $I_a(w_a)$ from above. By using the facts that the domain $\Omega_a$ is symmetric in $x$ and that the function $\chi_a$ is even in $x$ and by decomposing the integral $I_a(w_a)$ into three subdomains, one gets that
\be\label{Ianvn}\baa{rcl}
I_a(w_a) & \!\!=\!\! & \displaystyle\int_{(-1,1)^2}\frac{|\nabla(\phi_0(y)+t_0w(x,y))|^2}{2}\,dx\,dy-\int_{(-1,1)^2}(\phi_0(y)+t_0w(x,y))\,dx\,dy\vspace{3pt}\\
& & \displaystyle+2\int_{(1,a-1)\times(-1,1)}\frac{|\nabla\phi_0(y)|^2}{2}\,dx\,dy-2\int_{(1,a-1)\times(-1,1)}\phi_0(y)\,dx\,dy\vspace{3pt}\\
& & \displaystyle+2\int_{(a-1,a)\times(-1,1)}\frac{|\nabla(\phi_0(y)\chi_a(x))|^2}{2}\,dx\,dy-2\int_{(a-1,a)\times(-1,1)}\phi_0(y)\chi_a(x)\,dx\,dy\vspace{3pt}\\
& \!\!=\!\! & 2(a-2)J(\phi_0)+\beta,\eaa
\ee
where $\beta$ is a real number which does not depend on $a$ (it is indeed immediate to see by setting~$x=x'+a$ in the last two integrals of~(\ref{Ianvn}) that these quantities do not depend on $a$). Finally, it follows from~(\ref{Ianunvn}) and~(\ref{Ianvn}) that
\be\label{Iaua}
I_a(u_a)\le2(a-2)J(\phi_0)+\beta.
\ee\par
In the second step, we bound $I_a(u_a)$ from below. On the set $\Omega_a\backslash(-a,a)\times(-1,1)$, one simply uses the fact that
$$\int_{\Omega_a\backslash(-a,a)\times(-1,1)}\Big(\frac{|\nabla u_a|^2}{2}-u_a\Big)\ge-\int_{\Omega_a\backslash(-a,a)\times(-1,1)}\frac{5}{2}\ge-10\|\varphi\|_{L^{\infty}(-1,1)}$$
from~(\ref{uasup}) and from the definition~(\ref{defOmegaa}) of $\Omega_a$. Therefore,
\be\label{Iuaninf}\baa{rcl}
I_a(u_a) & \ge & \displaystyle\int_{(-a,a)\times(-1,1)}\!\!\Big(\frac{|\nabla u_a|^2}{2}-u_a\Big)-10\|\varphi\|_{L^{\infty}(-1,1)}\vspace{3pt}\\
& \ge & \displaystyle\int_{-a}^{a}\!\!J(u_a(x,\cdot))\,dx-10\|\varphi\|_{L^{\infty}(-1,1)},\eaa
\ee
where the functional $J$ has been defined in~(\ref{defJ}) and where we have used the fact that~$u_a(x,\cdot)$ belongs to~$H^1_0(-1,1)$ for all $x\in(-a,a)$. Remember that $\phi_0$ is the (unique) minimizer of~$J$. As a consequence,
\be\label{Juaphi}
J(u_a(x,\cdot))\ge J(\phi_0)\hbox{ for all }x\in(-a,a).
\ee
On the other hand, by definition of $x_a$ in~(\ref{defxa}) and by convexity and symmetry of $\omega_a$ with respect to both variables $x$ and $y$, it follows that $(x,0)\in\omega_a$ for all $x\in(-x_a,x_a)$, whence
$$u_a(x,0)>1>\phi_0(0)\ \hbox{ for all }x\in(-x_a,x_a).$$
Hence, there is a positive real number $\gamma>0$, independent of $a$, such that
$$\|u_a(x,\cdot)-\phi_0\|_{H^1(-1,1)}\ge\gamma>0\ \hbox{ for all }x\in(-x_a,x_a).$$
By definition of $\phi_0$ and from the coercivity of the functional $J$, one infers the existence of a positive constant $\delta>0$, independent of $a$, such that
$$J(u_a(x,\cdot))\ge J(\phi_0)+\delta\ \hbox{ for all }x\in(-x_a,x_a).$$
From~(\ref{Iuaninf}) and~(\ref{Juaphi}), one then gets that
\be\label{Iuaninf2}
I_a(u_a)\ge2\delta\,\min(x_a,a)+2aJ(\phi_0)-10\|\varphi\|_{L^{\infty}(-1,1)}.
\ee\par
Putting together~(\ref{Iaua}) and~(\ref{Iuaninf2}) with the inequality $x_a-\|\varphi\|_{L^{\infty}(-1,1)}\le\min(x_a,a)$ yields 
$$2\delta(x_a-\|\varphi\|_{L^{\infty}(-1,1)})+2aJ(\phi_0)-10\|\varphi\|_{L^{\infty}(-1,1)}\le2(a-2)J(\phi_0)+\beta,$$
where $\beta>0$ and $\delta>0$ are independent of $a$. Hence, there exists a constant $C_x>0$, independent of $a$, such that $0\le x_a<C_x$, that is~(\ref{taille1}). The proof of Lemma~\ref{lemtaille1} is thereby complete.\hfill$\Box$\break

The second lemma gives a bound from below of the ``vertical" size of the sets $\omega_a$, meaning that the sets~$\omega_a$ are not too thin.

\begin{lem}\label{lemtaille2}
There exists a constant $C_y>0$ such that 
\be\label{taille2}
0<C_y<\sup_{(x,y)\in\omega_a}|y|
\ee
for all $a\ge 1$ and for any minimizer $u_a$ of $I_a$ in $U_a$.
\end{lem}

\noindent{\bf{Proof.}} It is actually an immediate consequence of Lemma~\ref{lemtaille1} and of the constraint in the definition of the sets $U_a$. Consider any real number $a\ge 1$ and any minimizer $u_a$ of the functional $I_a$ in the set~$U_a$. Denote
\be\label{defya}
y_a=\sup_{(x,y)\in\omega_a}|y|.
\ee
There holds $y_a>0$ since $\omega_a$ is open and non-empty. Nevertheless, we want to get a lower bound that is independent of $a$. By Lemma~\ref{lemtaille1} and by definition of $x_a$ and $y_a$ in~(\ref{defxa}) and~(\ref{defya}), there holds
$$u_a\le 1\ \hbox{ in }\Omega_a\,\backslash\,(-C_x,C_x)\times(-y_a,y_a),$$
whence $g(u_a)=0$ in this set, using~(\ref{choiceg}). Therefore, since $u_a\in U_a$ and $g\le 1$ in $\R$, it follows that
$$1=\int_{\Omega_a}g(u_a)=\int_{\Omega_a\,\cap\,(-C_x,C_x)\times(-y_a,y_a)}g(u_a)\le4\,C_x\,y_a.$$
In other words, the conclusion~(\ref{taille2}) holds with $C_y$ such that $0<C_y<(4C_x)^{-1}$.\hfill$\Box$\break

{\it Step 4: comparison of $u_a(x,y)$ with $\phi_0(y)$ when $a$ is large.} In this step, we prove that the minimizers $u_a$ of $I_a$ in $U_a$ are close to the one-dimensional profile $\phi_0(y)=(1-y^2)/2$ far away from the origin and far away from the leftmost and rightmost points of $\overline{\Omega_a}$ in the direction $x$.

\begin{lem}\label{lemuaphi0}
For all $\epsilon>0$, there exist $A\ge 1$ and $M\in[0,A/2]$ such that
$$|u_a(x,y)-\phi_0(y)|=\Big|u_a(x,y)-\frac{1-y^2}{2}\Big|\le\epsilon\hbox{ in }\big([-a+M,-M]\cup[M,a-M]\big)\times[-1,1]\ (\subset\overline{\Omega_a}),$$
for all $a\ge A$ and for any minimizer $u_a$ of $I_a$ in $U_a$.
\end{lem}

\noindent{\bf{Proof.}} Assume that the conclusion does not hold for some $\epsilon>0$. Then there are some sequences~$(a_n)_{n\in\N}$ and $(x_n,y_n)_{n\in\N}$ of real numbers and points in $\R^2$ such that
\be\label{xnyn}
a_n\ge n,\ \ \frac{n}{2}\le|x_n|\le a_n-\frac{n}{2},\ \ |y_n|\le 1,\ \ |u_{a_n}(x_n,y_n)-\phi_0(y_n)|>\epsilon\ \hbox{ for all }n\in\N,
\ee
where $u_{a_n}$ is a minimizer of the functional $I_{a_n}$ in the set $U_{a_n}$. For each $n\in\N$, define
$$\overline{u}_n(x,y)=u_{a_n}(x+x_n,y)\ \hbox{ for all }(x,y)\in\overline{\Omega_{a_n}}-(x_n,0)=\big\{(x,y)\in\R^2;\ (x+x_n,y)\in\overline{\Omega_{a_n}}\big\}.$$
Each function $\overline{u}_n$ satisfies a semilinear elliptic equation of the type~(\ref{equa}) in $\overline{\Omega_{a_n}}-(x_n,0)$ with a nonlinearity~$f_{a_n}=1+\mu_{a_n}g'$ for some $\mu_{a_n}\in\R$. Lemma~\ref{lemtaille1} and~(\ref{vaua}) imply that
$$0<u_{a_n}(x,y)\le 1\ \hbox{ for all }n\in\N\hbox{ and }(x,y)\in\Omega_{a_n}\backslash(-C_x,C_x)\times(-1,1).$$
Hence, because of~(\ref{choiceg}) and~(\ref{xnyn}), for every fixed $C\ge 0$, there holds
$$0\le\overline{u}_n(x,y)\le 1\hbox{ and }\Delta\overline{u}_n(x,y)+1=0\hbox{ for all }(x,y)\in[-C,C]\times[-1,1],$$
for all $n$ large enough. From standard elliptic estimates up to the boundary, it follows that, up to extraction of a subsequence, the functions $\overline{u}_n$ converge in $C^2_{loc}(\R\times[-1,1])$ to a classical solution~$\overline{u}_{\infty}$ of
$$\left\{\baa{rclll}
\Delta\overline{u}_{\infty}+1 & = & 0 & \hbox{in} & \R\times[-1,1],\vspace{3pt}\\
0\ \le \ \overline{u}_{\infty} & \le & 1 & \hbox{in} & \R\times[-1,1],\vspace{3pt}\\
\overline{u}_{\infty} & = & 0 & \hbox{on} & \R\times\{\pm 1\}.\eaa\right.$$
Withtout loss of generality, one can also assume that $y_n\to y_{\infty}\in[-1,1]$ as $n\to+\infty$, whence
\be\label{overuinfty}
|\overline{u}_{\infty}(0,y_{\infty})-\phi_0(y_{\infty})|\ge\epsilon
\ee
from~(\ref{xnyn}).
\par
On the other hand, a standard Liouville-type result implies that $\overline{u}_{\infty}$ is necessarily identically equal to the one-dimensional profile $\phi_0(y)$ in $\R\times[-1,1]$. Indeed, the function $h(x,y)=\overline{u}_{\infty}(x,y)-\phi_0(y)$ is bounded and harmonic in $\R\times[-1,1]$, and it vanishes on $\R\times\{\pm 1\}$. The maximum principle implies that
$$|h(x,y)|\le\eta\,\cos\Big(\frac{\pi y}{4}\Big)\,\cosh\Big(\frac{\pi x}{4}\Big)$$
for all $(x,y)\in\R\times[-1,1]$ and for all $\eta>0$ (otherwise, the inequality would hold in $\R\times[-1,1]$ for some $\eta^*>0$, with equality at some point in $\R\times(-1,1)$, contradicting the strong maximum principle). Thus, since $\eta>0$ can be arbitrarily small, one gets that $h(x,y)=0$ for all $(x,y)\in\R\times[-1,1]$. In other words,
$$\overline{u}_{\infty}(x,y)=\phi_0(y)\ \hbox{ for all }(x,y)\in\R\times[-1,1].$$
This is in contradiction with~(\ref{overuinfty}) and the proof of Lemma~\ref{lemuaphi0} is thereby complete.\hfill$\Box$\break

{\it Step 5: the superlevel sets of the minimizers $u_a$ cannot be all convex when $a$ is large enough.} In this last step, we complete the proof of Theorem~\ref{th1}. Actually, Lemma~\ref{lemtaille2} and the one-dimensional convergence given in Lemma~\ref{lemuaphi0} will prevent any minimizer~$u_a$ of~$I_a$ in~$U_a$ from being quasiconcave when~$a$ is large enough.\par
Given $C_y>0$ as in Lemma~\ref{lemtaille2}, let $P$, $Q_a$ and $R_a$ be the points of $\R^2$ whose coordinates are given by
$$P=(0,C_y),\ Q_a=\Big(\frac{a}{4},\frac{C_y}{2}\Big)\hbox{ and }R_a=\Big(\frac{a}{2},0\Big)$$
for all $a\ge 1$, see the joint figure. From Lemma~\ref{lemtaille2} and the convexity and symmetry of $\omega_a$ with respect to $x$ and $y$, there holds $P\in\omega_a$, that is
$$u_a(P)>1$$
for any minimizer~$u_a$ of~$I_a$ in~$U_a$. On the other hand, the point $R_a$ belongs to $\Omega_a$ for all $a\ge 1$ by definition~(\ref{defOmegaa}) of~$\Omega_a$ and the point $Q_a$ is at the middle of the segment $[P,R_a]$ and is thus in $\Omega_a$ too by convexity of $\Omega_a$.\par
\begin{figure}\label{figure3}
\begin{center}
\includegraphics[scale=0.75]{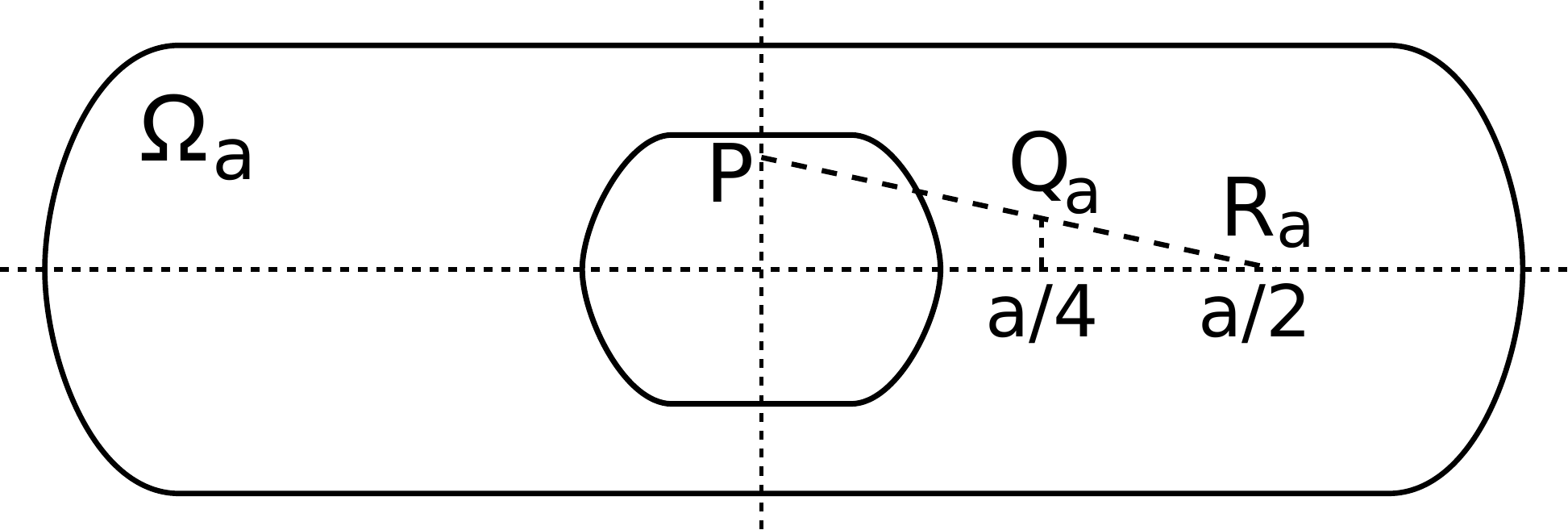}
\caption{The aligned points $P$, $Q_a$ and $R_a$}
\end{center}
\end{figure}
Furthermore, Lemma~\ref{lemuaphi0} implies that
$$u_a(Q_a)\longrightarrow\frac{1-(C_y/2)^2}{2}=\frac{1}{2}-\frac{C_y^2}{8}\ \hbox{ and }\ u_a(R_a)\longrightarrow\frac{1}{2}\ \hbox{ as }a\to+\infty,$$
for any minimizer~$u_a$ of~$I_a$ in~$U_a$. As a consequence, given any real number $\lambda$ such that
$$\frac{1}{2}-\frac{C_y^2}{8}<\lambda<\frac{1}{2},$$
one has
$$u_a(Q_a)<\lambda<u_a(R_a)<1<u_a(P)\ \hbox{ for all }a\hbox{ large enough}$$
and for any minimizer~$u_a$ of~$I_a$ in~$U_a$. Since the point $Q_a$ belongs to the segment $[P,R_a]$, it follows that, for $a$ large enough, the superlevel set
\be\label{upperlevel}
\big\{(x,y)\in\Omega_a;\ u_a(x,y)>\lambda\big\}
\ee
of any minimizer~$u_a$ of~$I_a$ in~$U_a$ is not convex, whence $u_a$ is not quasiconcave. The proof of Theorem~\ref{th1} is thereby complete.\hfill$\Box$

\begin{rem}{\rm By replacing $Q_a$ by $\widetilde{Q}_a=(\epsilon a,(1-2\epsilon)C_y)$ and by choosing $\epsilon\in(0,1/2)$ arbitrarily small, it follows from the above arguments that, given any real number $\lambda$ such that
$$\frac{1-C_y^2}{2}<\lambda<\frac{1}{2},$$
the superlevel set~$(\ref{upperlevel})$ of any minimizer~$u_a$ of~$I_a$ in~$U_a$ is not convex when $a$ is large enough.}
\end{rem}

\begin{rem}\label{remf'}{\rm In connection with Remark~$\ref{remstability}$ on the question of the stability of $u$, we focus here on the question of the sign of $f'_a$ on the range of $u_a$, for any minimizer $u_a$ of $I_a$ in $U_a$. First, we observe that $u_a\le 1$ and $f'_a(u_a)=0$ in $\Omega_a\,\backslash\,(-C_x,C_x)\times(-1,1)$, which is a large set when $a$ is large. Furthermore, the set
$$E^+_a=\big\{(x,y)\in\Omega_a;\ f'_a(u_a(x,y))>0\big\}$$
is never empty. Indeed, if $E^+_a$ were empty, then $g''=f'_a/\mu_a$ would be nonpositive on the range of~$u_a$, that is on the interval $[0,\max_{\overline{\Omega_a}}u_a]$. Due to~$(\ref{choiceg})$, that would mean that $g$ is actually constant equal to~$0$ on this interval $[0,\max_{\overline{\Omega_a}}u_a]$, whence $f_a(u_a)=1+\mu_ag'(u_a)=1$ in $\Omega_a$. That would imply that~$u_a=v_a$ in~$\Omega_a$, which is not the case. Thus, the open set~$E^+_a$ cannot be the empty set. On the other hand, the set
$$E^-_a=\big\{(x,y)\in\Omega_a;\ f'_a(u_a(x,y))<0\big\}$$
is not empty in general. Indeed, let for instance $\theta$ be the function defined in $(1,2)$ by
$$\theta(s)=e^{-\frac{1}{s-1}-\frac{1}{2-s}}\Big(\sin\Big(\frac{1}{(s-1)^2}\Big)+1\Big),\ \ s\in(1,2)$$
and let the function $g$ be defined by
$$g(s)=\left\{\baa{ll} 0 & \hbox{if }s\le 1,\vspace{3pt}\\
\kappa\displaystyle\int_1^s\theta(t)\,dt & \hbox{if }1<s<2,\vspace{3pt}\\
1 & \hbox{if }s\ge 2,\eaa\right.$$
where the constant $\kappa>0$ is chosen so that $g$ is continuous at $s=2$. The function $g$ is then of class~$C^{\infty}(\R)$ and it satisfies~$(\ref{choiceg})$. But $g''$ has infinitely many sign changes in any right neighborhood of $1$. For this choice of $g$ and for any minimizer $u_a$ of $I_a$ in $U_a$ with $a\ge 1$, since~$0=\min_{\overline{\Omega_a}}u_a<1<\max_{\overline{\Omega_a}}u_a$ and $f'_a(u_a)=\mu_ag''(u_a)$, it follows that the set $E^-_a$ is not empty.}
\end{rem}

%%%%%%%%%%%%%%%%%%%%%%%%%%%%%%%%%%%%%%%%%
%%%%%%%%%%%%%%%%%%%%%%%%%%%%%%%%%%%%%%%%%

\section{Counterexamples in convex rings}\label{secrings}

In this section, we consider problems of the type~(\ref{pbring}) set in convex rings $\Omega=\Omega_1\backslash\overline{\Omega_2}$. The examples of non-convexity of some superlevel sets of the solutions of~(\ref{pbring}) stated in Theorem~\ref{th2} can be viewed as a particular case of a more general statement. Namely, we shall construct counterexamples for the convexity of the superlevel sets of the solutions of heterogeneous non-symmetric semilinear elliptic equations of the type
\be\label{pbringbis}\left\{\baa{rcll}
\nabla\cdot(A(x)\nabla u)+b(x)\cdot\nabla u+f(x,u) & \!\!=\!\! & 0 & \hbox{in }\Omega,\vspace{3pt}\\
u & \!\!=\!\! & 0 & \hbox{on }\partial\Omega_1,\vspace{3pt}\\
u & \!\!=\!\! & M & \hbox{on }\partial\Omega_2,\vspace{3pt}\\
u & \!\!>\!\! & 0 & \hbox{in }\Omega,\eaa\right.
\ee
where $M>0$ is a positive real number and $\Omega=\Omega_1\backslash\overline{\Omega_2}$ is a convex ring. The convex domain~$\Omega_1$ is given and its boundary is smooth, in the sense that it is of class $C^{2,\alpha}$ with~$\alpha>0$. The convex interior domain $\Omega_2$ such that $\overline{\Omega_2}\subset\Omega_1$ shall be constructed later, in the proof of Theorem~\ref{th2bis} below. The coefficients $A$ and $b$ are given in~$\overline{\Omega_1}$ and $f$ in $\overline{\Omega_1}\times[0,+\infty)$. More precisely, the matrix field~$A:x\mapsto A(x)=(a_{ij}(x))_{1\le i,j\le N}$ is a symmetric matrix field of class~$C^{1,\alpha}(\overline{\Omega_1})$ such that
$$\exists\,\beta>0,\ \forall\,x\in\overline{\Omega_1},\ \forall\,\xi=(\xi_i)_{1\le i\le N}\in\R^N,\ \sum_{1\le i,j\le N}a_{ij}(x)\,\xi_i\,\xi_j\ge\beta\,|\xi|^2,$$
where $|\xi|^2=\xi_1^2+\cdots+\xi_N^2$. The vector field $b:x\mapsto b(x)=(b_i(x))_{1\le i\le N}$ is of class $C^{0,\alpha}(\overline{\Omega_1})$ and the function $f:\overline{\Omega_1}\times[0,+\infty),\ (x,s)\mapsto f(x,s)$ is of class $C^{0,\alpha}(\overline{\Omega_1})$ with respect to~$x$ locally in~$s$ and locally Lipschitz-continuous with respect to $s$ uniformly in $x$. Furthermore, we assume that
\be\label{hypfbis}\left\{\baa{l}
f\hbox{ is bounded from above, that is }\displaystyle{\mathop{\sup}_{(x,s)\in\overline{\Omega_1}\times[0,+\infty)}}\,f(x,s)<+\infty,\vspace{3pt}\\
s\mapsto\displaystyle{\frac{f(x,s)}{s}}\hbox{ is nonincreasing over }(0,+\infty)\hbox{ for all }x\in\Omega_1,\vspace{3pt}\\
s\mapsto\displaystyle{\frac{f(\tilde{x},s)}{s}}\hbox{ is decreasing over }(0,+\infty)\hbox{ for some }\tilde{x}\in\Omega_1,\vspace{3pt}\\
\hbox{either }\displaystyle{\mathop{\max}_{\overline{\Omega_1}}}\,f(\cdot,0)>0,\ \hbox{ or }\ f(\cdot,0)=0\hbox{ in }\overline{\Omega_1}\hbox{ and }\lambda_1(-\L,\Omega_1)<0,\eaa\right.
\ee
where
$$\left\{\baa{rcl}
\L & = & \nabla\cdot(A(x)\nabla)+b(x)\cdot\nabla+\zeta(x),\vspace{3pt}\\
\zeta(x) & = & \displaystyle{\mathop{\lim}_{s\to0^+}}\,\displaystyle{\frac{f(x,s)}{s}}\eaa\right.$$
and $\lambda_1(-\L,\Omega_1)$ denotes the principal eigenvalue of the operator $-\L$ in $\Omega_1$ with Dirichlet boundary condition on $\partial\Omega_1$. In the case $f(\cdot,0)=0$ in $\overline{\Omega_1}$, we assume moreover that $\zeta$ is H\"older-continuous in~$\overline{\Omega_1}$, whence the limit $\zeta=\lim_{s\to0^+}f(\cdot,s)/s$ is uniform in $\overline{\Omega_1}$ from Dini's theorem. In this case, the principal eigenvalue $\lambda_1(-\L,\Omega_1)$ of the operator $-\L$ with Dirichlet boundary condition on~$\partial\Omega_1$ is a real number which is characterized by the existence and uniqueness (up to multiplication) of a classical eigenfunction $\varphi$ solving
\be\label{varphi}\left\{\baa{rcll}
-\L\varphi & = & \lambda_1(-\L,\Omega_1)\,\varphi & \hbox{in }\overline{\Omega_1},\vspace{3pt}\\
\varphi & = & 0 & \hbox{on }\partial\Omega_1,\vspace{3pt}\\
\varphi & > & 0 & \hbox{in }\Omega_1,\eaa\right.
\ee
see~\cite{bnv}.\footnote{When~$\max_{\overline{\Omega_1}}f(\cdot,0)>0$, if we define $\L_s=\nabla\cdot(A(x)\nabla)+b(x)\cdot\nabla+f(x,s)/s$ for all $s>0$, the map~$s\mapsto\lambda_1(-\L_s,\Omega_1)$ is nondecreasing on $(0,+\infty)$ and one can then set $\lambda_1(-\L,\Omega_1)=\lim_{s\to0^+}\lambda_1(-\L_s,\Omega_1)$. On the other hand, $f(\cdot,s)/s\to+\infty$ as $s\to0^+$ at least uniformly in a subdomain of $\Omega_1$. Since $\lambda_1(-\L_s,\cdot)$ is nonincreasing with respect to the inclusion of domains for each $s>0$ (see~\cite{bnv}), it follows that $\lambda_1(-\L,\Omega_1)$, as defined as above, is equal to~$-\infty$.}\par
Theorem~\ref{th2} can then be viewed as a particular case of the following result.

\begin{theo}\label{th2bis}
Let $N$ be any integer such that $N\ge 2$, let $\Omega_1$ be any smooth bounded convex domain of~$\R^N$, let $A$ and $b$ be as above and let $f$ be any function satisfying~$(\ref{hypfbis})$ with $\zeta$ being H\"older-continuous in $\overline{\Omega_1}$ in case $f(\cdot,0)=0$ in $\overline{\Omega_1}$. Then there exists a constant $M_0>0$ such that, for all $M\ge M_0$, there are some smooth convex rings $\Omega=\Omega_1\backslash\overline{\Omega_2}$ for which problem~$(\ref{pbringbis})$ has a unique solution~$u$ and this solution~$u$ is not quasiconcave.
\end{theo}

\noindent{\bf{Proof.}} Let $N$, $\Omega_1$, $A$, $b$ and $f$ be as in the statement. The strategy of the proof consists in the following steps: we first construct and prove the uniqueness of a solution~$v$ of the boundary value problem
\be\label{pbOmega1}\left\{\baa{rcll}
\nabla\cdot(A(x)\nabla v)+b(x)\cdot\nabla v+f(x,v) & \!\!=\!\! & 0 & \hbox{in }\Omega_1,\vspace{3pt}\\
v & \!\!=\!\! & 0 & \hbox{on }\partial\Omega_1,\vspace{3pt}\\
v & \!\!>\!\! & 0 & \hbox{in }\Omega_1;\eaa\right.
\ee
next, for $M_0=\max_{\overline{\Omega_1}}v$, for any point $x_0\in\Omega_1$ such that $v(x_0)<\max_{\overline{\Omega_1}}v$ and for $\Omega_2$ being a smooth convex domain included in the Euclidean ball $B(x_0,\epsilon)$ of center $x_0$ and radius $\epsilon>0$ small enough, we prove the existence and uniqueness of a solution $u$ of~(\ref{pbringbis}) with $\Omega=\Omega_1\backslash\overline{\Omega_2}$ and~$M\ge M_0$; by uniqueness of $v$, this solution~$u$ shall be close to $v$ locally in $\overline{\Omega_1}\backslash\{x_0\}$ for~$\epsilon>0$ small enough, from standard elliptic estimates and a priori bounds; the conclusion, that is $u$ has some non-convex superlevel sets for $\epsilon>0$ small enough, will then follow from the choice of~$x_0$ and the fact that $M\ge\max_{\overline{\Omega_1}}v$.\hfill\break

{\it{Step 1: problem~$(\ref{pbOmega1})$ in~$\Omega_1$}}. Let us first prove the existence and uniqueness of a solution~$v$ of~(\ref{pbOmega1}) in~$\Omega_1$. The proof draws its inspiration from~\cite{b,bhr1,bhr2}, where $f$ is usually assumed to be nonpositive for $s$ large enough (instead of being globally bounded from above). We adapt the method with the weaker assumptions~(\ref{hypfbis}).\par
Let~$\psi$ be the unique~$C^{2,\alpha}(\overline{\Omega_1})$ solution of the boundary value problem
\be\label{eqpsi}\left\{\baa{rcll}
\nabla\cdot(A(x)\nabla\psi)+b(x)\cdot\nabla\psi & \!\!=\!\! & -1 & \hbox{in }\Omega_1,\vspace{3pt}\\
\psi & \!\!=\!\! & 0 & \hbox{on }\partial\Omega_1.\eaa\right.
\ee
The function $\psi$ is such that $\psi>0$ in $\Omega_1$ from the strong maximum principle. Let $C$ be a positive real number such that
\be\label{defC}
f(x,s)\le C\hbox{ for all }(x,s)\in\overline{\Omega_1}\times[0,+\infty).
\ee
It follows that the function $C\psi$ is a supersolution of the equation~(\ref{pbOmega1}) in~$\Omega_1$, in the sense that
$$\baa{rcl}
\nabla\cdot(A(x)\nabla(C\psi))+b(x)\cdot\nabla(C\psi)+f(x,C\psi) & \le &\nabla\cdot(A(x)\nabla(C\psi))+b(x)\cdot\nabla(C\psi)+C\vspace{3pt}\\
& = & 0\ \hbox{ in }\Omega_1.\eaa$$
Furthermore, Hopf's lemma implies that $\partial\psi/\partial\nu<0$ on $\partial\Omega_1$, where~$\nu$ denotes the outward unit normal on~$\partial\Omega_1$.\par
In order to construct a subsolution of~(\ref{pbOmega1}), we first consider the case when $f(\cdot,0)=0$ in~$\overline{\Omega_1}$ and $\lambda_1(-\L,\Omega_1)<0$. Let~$\varphi$ be a classical solution of the eigenvalue problem~(\ref{varphi}) in~$\overline{\Omega_1}$. Since the convergence $f(\cdot,s)/s\to\zeta$ as $s\to0^+$ is uniform in $\overline{\Omega_1}$ and since $\lambda_1(-\L,\Omega_1)<0$, it follows that there exists $\delta_0>0$ such that, for all $\delta\in(0,\delta_0)$,
\be\label{subvarphi}\left\{\baa{rcll}
\nabla\cdot(A(x)\nabla(\delta\varphi))+b(x)\cdot\nabla(\delta\varphi)+f(x,\delta\varphi) & \!\!\ge\!\! & 0 & \hbox{in }\Omega_1,\vspace{3pt}\\
\delta\varphi & \!\!=\!\! & 0 & \hbox{on }\partial\Omega_1,\vspace{3pt}\\
\delta\varphi & \!\!>\!\! & 0 & \hbox{in }\Omega_1.\eaa\right.
\ee
In other words, $\delta\varphi$ is a subsolution of problem~(\ref{pbOmega1}) in~$\Omega_1$ for $\delta>0$ small enough. Furthermore, since $\psi>0$ in $\Omega_1$, $\partial\psi/\partial\nu<0$ on $\partial\Omega_1$, $\psi$ and $\varphi$ both vanish on $\partial\Omega_1$ and $\varphi$ is (at least) of class~$C^1(\overline{\Omega_1})$, there holds
$$\delta\varphi\le C\psi\hbox{ in }\overline{\Omega_1}$$
for $\delta>0$ small enough. For some given small enough $\delta>0$, the monotone iteration method yields the existence of a solution~$v$ of~(\ref{pbOmega1}) such that
$$\delta\varphi\le v\le C\psi\hbox{ in }\overline{\Omega_1}.$$\par
Consider now the case when $\max_{\overline{\Omega_1}}f(\cdot,0)>0$. Let $B$ be a non-empty open Euclidean ball such that $\overline{B}\subset\Omega_1$ and $\min_{\overline{B}}f(\cdot,0)>0$. Let $\phi$ be any $C^2(\overline{B})$ function such that $\phi>0$ in~$B$ and~$\phi=0$ on~$\partial B$. There exists then $\tilde{\delta}_0>0$ such that, for all $\delta\in(0,\tilde{\delta}_0)$,
\be\label{subphi}\left\{\baa{rcll}
\nabla\cdot(A(x)\nabla(\delta\phi))+b(x)\cdot\nabla(\delta\phi)+f(x,\delta\phi) & \!\!\ge\!\! & 0 & \hbox{in }B,\vspace{3pt}\\
\delta\phi & \!\!=\!\! & 0 & \hbox{on }\partial B,\vspace{3pt}\\
\delta\phi & \!\!>\!\! & 0 & \hbox{in }B.\eaa\right.
\ee
On the other hand, there holds
\be\label{fge0}
f(x,0)\ge 0\hbox{ for all }x\in\overline{\Omega_1}
\ee
because $f(x,s)/s$ is nonincreasing in $s\in(0,+\infty)$ for all $x\in\overline{\Omega_1}$. Hence, the function $\delta\phi$ extended by~$0$ in~$\overline{\Omega_1}\backslash\overline{B}$ is a subsolution of problem~(\ref{pbOmega1}) in~$\Omega_1$, for $\delta>0$ small enough. Furthermore,
$$\delta\phi\le C\psi\hbox{ in }\overline{B}$$
for $\delta>0$ small enough. As above, one then gets the existence of a solution~$v$ of~(\ref{pbOmega1}) such that
$$\left\{\baa{l}
v\le C\psi\hbox{ in }\overline{\Omega_1},\vspace{3pt}\\
v\ge\delta\phi\hbox{ in }\overline{B}\hbox{ and }v\ge0\hbox{ in }\overline{\Omega_1}\backslash\overline{B}\eaa\right.$$
for some given small enough $\delta>0$. Notice in particular that $v>0$ in~$\Omega_1$ from~(\ref{fge0}) and the strong maximum principle.\par
Lastly, let us prove the uniqueness of the solution~$v$ of~(\ref{pbOmega1}). Let $w$ be another solution of~(\ref{pbOmega1}). Since~$v$ and $w$ are at least of class~$C^2(\overline{\Omega_1})$ and the constant~$0$ is always a subsolution of problem~(\ref{pbOmega1}) (because $f(\cdot,0)\ge0$ in~$\overline{\Omega_1}$), Hopf's lemma implies that $\partial v/\partial\nu<0$ and~$\partial w/\partial\nu<0$ on $\partial\Omega_1$. It follows that there exists a constant $\tau\ge 1$ such that
$$\tau^{-1}w\le v\le \tau w\hbox{ in }\overline{\Omega_1}.$$
Let $t^*\in[\tau^{-1},\tau]$ be defined as
$$t^*=\min\big\{t>0,\ v\le t\,w\hbox{ in }\overline{\Omega_1}\big\}.$$
Assume that $t^*>1$. Since $f(x,s)/s$ is nonincreasing with respect to $s\in(0,+\infty)$ for all $x\in\overline{\Omega_1}$ and decreasing for at least a point $\tilde{x}$ in $\Omega_1$, it follows that
\be\label{strictsuper}\baa{l}
\nabla\cdot(A(x)\nabla(t^*w))+b(x)\cdot\nabla(t^*w)+f(x,t^*w)\vspace{3pt}\\
\qquad\qquad\qquad\qquad\qquad\qquad\le,\not\equiv t^*\big(\nabla\cdot(A(x)\nabla w)+b(x)\cdot\nabla w+f(x,w)\big)=0\eaa
\ee
in $\Omega_1$, while $v\le t^*w$ in $\Omega_1$. One infers from the strong maximum principle that either $v<t^*w$ in~$\Omega_1$ or $v=t^*w$ in $\overline{\Omega_1}$. The first case is impossible since it would then imply that $v\le(t^*-\epsilon)w$ in $\overline{\Omega_1}$ for all $\epsilon>0$ small enough, using again Hopf's lemma, and it would contradict the definition of $t^*$. Thus,~$v=t^*w$ in $\overline{\Omega_1}$, which is also impossible since the inequality~(\ref{strictsuper}) is not an equality everywhere. As a consequence, $t^*\le 1$, whence
$$v\le w\hbox{ in }\overline{\Omega_1}.$$
Reversing the roles of $v$ and $w$ leads to the conclusion $v=w$ in $\overline{\Omega_1}$.\hfill\break

{\it{Step 2: Problem~$(\ref{pbringbis})$ in suitable convex rings $\Omega=\Omega_1\backslash\overline{\Omega_2}$.}} Set $M_0=\max_{\overline{\Omega_1}}v>0$ and pick any constant $M$ such that
\be\label{choiceM}
M\ge M_0=\max_{\overline{\Omega_1}}\,v.
\ee
Let us prove the existence of convex smooth domains $\Omega_2$ such that $\overline{\Omega_2}\subset\Omega_1$ and for which problem~(\ref{pbringbis}) has a unique solution~$u$ in $\Omega=\Omega_1\backslash\overline{\Omega_2}$ and this solution has some non-convex superlevel sets. To do so, pick any point $x_0\in\Omega_1$ such that
\be\label{choicex0}
v(x_0)<\max_{\overline{\Omega_1}}\,v,
\ee
let $\omega_2$ be any (smooth, that is of class $C^{2,\alpha}$) convex domain of $\R^N$ and consider convex rings of the type
$$\Omega^{\epsilon}=\Omega_1\backslash\overline{\Omega_2^{\epsilon}},\ \hbox{ with }\Omega_2^{\epsilon}=x_0+\epsilon\,\omega_2$$
for $\epsilon>0$ small enough: namely, there is $\epsilon^*>0$ such that $\overline{\Omega_2^{\epsilon}}\subset\Omega_1$ and $\Omega^{\epsilon}$ is then a  convex ring for all $\epsilon\in(0,\epsilon^*)$. Without loss of generality, one can also assume that there is a fixed real number~$r>0$ such that
$$\overline{\Omega_2^{\epsilon}}\subset\overline{B(x_0,r)}\subset\Omega_1\hbox{ for all }\epsilon\in(0,\epsilon^*),$$
where $B(x_0,r)$ denotes the open Euclidean ball of center $x_0$ and radius $r>0$.\par
Remember that $\psi$ is defined by~(\ref{eqpsi}) in $\overline{\Omega_1}$. Since $\psi$ is continuous and positive in $\Omega_1$, there holds~$D\psi\to+\infty$ locally uniformly in $\Omega_1$ as $D\to+\infty$ (in particular, $\min_{\overline{B(x_0,r)}}D\psi\to+\infty$ as~$D\to+\infty$). On the other hand, for all $D\ge C$, one has
$$\nabla\cdot(A(x)\nabla(D\psi))+b(x)\cdot\nabla(D\psi)+f(x,D\psi)\le -D+C\le 0\hbox{ in }\Omega_1$$
because of~(\ref{eqpsi}) and~(\ref{defC}). Therefore, there is a positive constant $D\ge C$ such that the function~$D\psi$ is a supersolution of problem~(\ref{pbringbis}) in the convex ring $\Omega^{\epsilon}$ for all $\epsilon\in(0,\epsilon^*)$ (in particular, one has $D\psi\ge M$ on $\partial\Omega^{\epsilon}_2$).\par
On the other hand, when $f(\cdot,0)=0$ in $\overline{\Omega_1}$, let $\varphi$ solve the eigenvalue problem~(\ref{varphi}) in~$\Omega_1$ and let $\delta>0$ be small enough so that
$$\delta\|\varphi\|_{L^{\infty}(\Omega_1)}<M,\ \ \delta\varphi\le D\psi\hbox{ in }\overline{\Omega_1}$$
and $\delta\varphi$ be a subsolution of~(\ref{pbOmega1}) in~$\Omega_1$, that is $\delta\varphi$ satisfies~(\ref{subvarphi}). Choosing such a $\delta>0$ is possible since $\psi>0$ in $\Omega_1$, $\partial\psi/\partial\nu<0$ on $\partial\Omega_1$ and $\varphi$ is (at least) of class $C^1(\overline{\Omega_1})$. The function~$\delta\varphi$ is then a subsolution of problem~(\ref{pbringbis}) in $\Omega^{\epsilon}$ for all $\epsilon\in(0,\epsilon^*)$. When $\max_{\overline{\Omega_1}}f(\cdot,0)>0$, let~$B$ and~$\phi$ be as in Step~1 and let $\delta>0$ small enough so that
$$\delta\|\phi\|_{L^{\infty}(B)}<M,\ \ \delta\phi\le D\psi\hbox{ in }\overline{B}$$
and $\delta\phi$ (extended by $0$ in $\overline{\Omega_1}\backslash\overline{B}$) be a subsolution of~(\ref{pbOmega1}) in~$\Omega_1$, that is $\delta\phi$ satisfies~(\ref{subphi}). The function $\delta\phi$ (extended by $0$ in $\overline{\Omega_1}\backslash\overline{B}$) is then a subsolution of problem~(\ref{pbringbis}) in $\Omega^{\epsilon}$ for all~$\epsilon\in(0,\epsilon^*)$.\par
In both cases $f(\cdot,0)=0$ in $\overline{\Omega_1}$ and $\max_{\overline{\Omega_1}}f(\cdot,0)>0$, for every $\epsilon\in(0,\epsilon^*)$, there exists a solution $u^{\epsilon}$ of~(\ref{pbringbis}) in $\Omega^{\epsilon}$ such that
\be\label{boundsueps}\left\{\baa{lcl}
f(\cdot,0)=0\hbox{ in }\overline{\Omega_1} & \Longrightarrow & \delta\varphi\le u^{\epsilon}\le D\psi\hbox{ in }\overline{\Omega^{\epsilon}},\vspace{3pt}\\
\displaystyle{\mathop{\max}_{\overline{\Omega_1}}}\,f(\cdot,0)>0 & \Longrightarrow & \delta\phi\le u^{\epsilon}\le D\psi\hbox{ in }\overline{\Omega^{\epsilon}}\cap\overline{B}\hbox{ and }0\le u^{\epsilon}\le D\psi\hbox{ in }\overline{\Omega^{\epsilon}}\backslash\overline{B}.\eaa\right.
\ee
In particular, since $u^{\epsilon}$ is nonnegative by construction and not identically equal to $0$ in $\Omega^{\epsilon}$ (because, for instance,~$u^{\epsilon}=M>0$ on $\partial\Omega_2^{\epsilon}$) and since $0$ is always a subsolution of problem~(\ref{pbringbis}) (because~$f(\cdot,0)\ge 0$ in $\overline{\Omega_1}$), the strong maximum principle yields $u^{\epsilon}>0$ in $\Omega^{\epsilon}$. Observe now that, if $v^{\epsilon}$ is another solution of~(\ref{pbringbis}) in~$\Omega^{\epsilon}$, then the equality~$u^{\epsilon}=t\,v^{\epsilon}$ in $\Omega^{\epsilon}$ for some~$t>0$ with $t\neq 1$ is impossible due to the boundary condition on $\partial\Omega_2^{\epsilon}$. Therefore, by using the same method as in Step~1, whether $\tilde{x}$ be in $\Omega^{\epsilon}_2$ or not, it follows that the solution $u^{\epsilon}$ of~(\ref{pbringbis}) in~$\Omega^{\epsilon}$ is unique.\hfill\break

{\it{Step 3: non-convexity of some superlevel sets of $u^{\epsilon}$ in $\Omega^{\epsilon}$ for $\epsilon>0$ small enough.}} We now complete the proof of Theorem~\ref{th2bis}. We first claim that
\be\label{claim}
u^{\epsilon}\to v\hbox{ in }C^2_{loc}(\overline{\Omega_1}\backslash\{x_0\})\hbox{ as }\epsilon\to0^+,
\ee
where $v$ denotes the unique solution of~(\ref{pbOmega1}), given in Step~1.\par
To prove this claim, let $(\epsilon_n)_{n\in\N}$ be any sequence of real numbers in $(0,\epsilon^*)$ such that $\epsilon_n\to0^+$ as~$n\to+\infty$. The sequence $(\|u^{\epsilon_n}\|_{L^{\infty}(\Omega^{\epsilon_n})})_{n\in\N}$ is bounded from~(\ref{boundsueps}) (remember that the constant~$D$ is independent of $\epsilon$). For any compact subset~$K\subset\overline{\Omega_1}\backslash\{x_0\}$, the sequence $(u^{\epsilon_n})_{n\ge n_0}$ is then bounded in $C^{2,\alpha}(K)$ for $n_0$ large enough, from standard elliptic estimates and from the definition of $\Omega^{\epsilon}$. Up to extraction of a subsequence, the functions $u^{\epsilon_n}$ converge as $n\to+\infty$ in~$C^2_{loc}(\overline{\Omega_1}\backslash\{x_0\})$ to a $C^2(\overline{\Omega_1}\backslash\{x_0\})$ solution $u^0$ of
$$\left\{\baa{rcll}
\nabla\cdot(A(x)\nabla u^0)+b(x)\cdot\nabla u^0+f(x,u^0) & \!\!=\!\! & 0 & \hbox{in }\Omega_1\backslash\{x_0\},\vspace{3pt}\\
u^0 & \!\!=\!\! & 0 & \hbox{on }\partial\Omega_1\eaa\right.$$
such that
$$0\le\delta\varphi\le u^0\le D\psi\hbox{ in }\overline{\Omega_1}\backslash\{x_0\}$$
when $f(\cdot,0)=0$ in $\overline{\Omega_1}$, resp.
$$0\le\delta\phi\le u^0\le D\psi\hbox{ in }\overline{B}\backslash\{x_0\}\hbox{ and }0\le u^0\le D\psi\hbox{ in }(\overline{\Omega_1}\backslash\overline{B})\backslash\{x_0\}$$
when $\max_{\overline{\Omega_1}}f(\cdot,0)>0$ (remember that the positive constants $\delta$ and $D$ are independent of $\epsilon$). Since the set~$\{x_0\}$ is removable~\cite{ar,s}, it follows that $u_0$ can be extended to a $C^2(\overline{\Omega_1})$ solution of~(\ref{pbOmega1}). In particular, notice that the positivity of $u^0$ in $\Omega_1$ follows from the strong maximum principle and the lower bound $u^0\ge\delta\varphi$ in $\overline{\Omega_1}$ when $f(\cdot,0)=0$ in $\overline{\Omega_1}$, resp. $u^0\ge\delta\phi$ in $\overline{B}$ when~$\max_{\overline{\Omega_1}}f(\cdot,0)>0$. From Step~1 and the uniqueness of the solution of~(\ref{pbOmega1}), one gets that
$$u^0=v\hbox{ in }\overline{\Omega_1},$$
whence $u^{\epsilon_n}\to v$ in $C^2_{loc}(\overline{\Omega_1}\backslash\{x_0\})$ as $n\to+\infty$. Since the limit does not depend on the sequence~$(\epsilon_n)_{n\in\N}$, the claim~(\ref{claim}) follows.\par
To get the conclusion of Theorem~\ref{th2bis}, it is then sufficient to prove that the solutions $u^{\epsilon}$ of~(\ref{pbringbis}) in~$\Omega^{\epsilon}$ have some non-convex superlevel sets, at least for $\epsilon>0$ small enough. Assume by contradiction that this conclusion does not hold, that is there is a sequence~$(\epsilon_n)_{n\in\N}$ in $(0,\epsilon^*)$ such that $\epsilon_n\to0^+$ as~$n\to+\infty$ and, for each $n\in\N$, the superlevel sets of the function $u^{\epsilon_n}$ are all convex. For each~$n\in\N$, extend the function $u^{\epsilon_n}$ by $M$ in $\Omega_2^{\epsilon_n}=x_0+\epsilon_n\omega_2$, and still call $u^{\epsilon_n}$ this extension, now defined in $\overline{\Omega_1}$. Fix a point $y\in\Omega_1$ such that
\be\label{defy}
v(y)=\max_{\overline{\Omega_1}}\,v>0
\ee
and let $(x_n)_{n\in\N}$ be any sequence of points in $\overline{\Omega_1}$ such that $x_n\in\overline{\Omega_2^{\epsilon_n}}$ for all $n\in\N$. Lastly, let~$\eta>0$ be an arbitrary positive real number. Since $y\neq x_0$ from~(\ref{choicex0}) and~(\ref{defy}), the convergence~(\ref{claim}) implies in particular that $u^{\epsilon_n}(y)\to v(y)$ as $n\to+\infty$. Therefore, there is $n_0\in\N$ such that
$$u^{\epsilon_n}(y)\ge v(y)-\eta\hbox{ for all }n\ge n_0.$$
On the other hand, $u^{\epsilon_n}(x_n)=M\ge M_0=\max_{\overline{\Omega_1}}v=v(y)$ from~(\ref{choiceM}) and~(\ref{defy}). Since the superlevel sets of $u^{\epsilon_n}$ in $\overline{\Omega_1}$ are assumed to be all convex, it follows that
$$u^{\epsilon_n}(x)\ge v(y)-\eta\hbox{ for all }x\in[x_n,y]\hbox{ and for all }n\ge n_0,$$
where $[x_n,y]$ denotes the segment between $x_n$ and $y$, that is
$$[x_n,y]=\{tx_n+(1-t)y,\ 0\le t\le1\}.$$
From~(\ref{claim}) and the fact that $x_n\to x_0$ as $n\to+\infty$, one infers that $v(x)\ge v(y)-\eta$ for all~$x\in(x_0,y]$, and then also at the point $x=x_0$ be continuity of $v$. Since $\eta>0$ is arbitrary, it follows that
$$v(x_0)\ge v(y)=\max_{\overline{\Omega_1}}\,v,$$
which is ruled out due to the choice of $x_0$ in~(\ref{choicex0}). One has then reached a contradiction and the proof of Theorem~\ref{th2bis} is thereby complete.\hfill$\Box$

\begin{rem}\label{remM}{\rm When, in addition to~$(\ref{hypfbis})$, the function $f$ is nonpositive for large $s$ uniformly in~$x$, that is there exists a constant $\mu>0$ such that
$$f(x,s)\le 0\hbox{ for all }(x,s)\in\overline{\Omega_1}\times[\mu,+\infty),$$
then any constant $M$ such that $M\ge\mu$ is a supersolution of problem~$(\ref{pbOmega1})$ in $\Omega_1$ and~$(\ref{pbringbis})$ in~$\Omega^{\epsilon}$ (for $\epsilon>0$ small enough). Therefore, from the strong maximum principle, for any constant~$M\ge\mu$, the unique solution $v$ of~$(\ref{pbOmega1})$ in $\Omega_1$ and the unique solution $u^{\epsilon}$ of~$(\ref{pbringbis})$ in $\Omega^{\epsilon}$ satisfy $v<M$ in~$\Omega_1$ and $u^{\epsilon}<M$ in $\Omega^{\epsilon}$ for $\epsilon>0$ small enough.}
\end{rem}

%%%%%%%%%%%%%%%%%%%%%%%%%%%%%%%%%%%%%%%%%
%%%%%%%%%%%%%%%%%%%%%%%%%%%%%%%%%%%%%%%%%

\end{document}